\documentclass{amsart}
\usepackage{amsxtra}
\theoremstyle{plain}

\newcommand{\nn}{\hfill\nonumber}
\newtheorem{theorem}{Theorem}[section]
\newtheorem{lemma}[theorem]{Lemma}
\newtheorem{definition-theorem}[theorem]{Definition-Theorem}
\newtheorem{proposition}[theorem]{Proposition}
\newtheorem{corollary}[theorem]{Corollary}
\newtheorem{definition}[theorem]{Definition}
\newtheorem{example}[theorem]{Example}
\newtheorem{remark}[theorem]{Remark}
\newcommand \bth[1] { \begin{theorem}\label{t#1} }
\newcommand \ble[1] { \begin{lemma}\label{l#1} }

\newcommand \bpr[1] { \begin{proposition}\label{p#1} }
\newcommand \bco[1] { \begin{corollary}\label{c#1} }
\newcommand \bde[1] { \begin{definition}\label{d#1}\rm }
\newcommand \bex[1] { \begin{example}\label{e#1}\rm }
\newcommand \bre[1] { \begin{remark}\label{r#1}\rm }

\newcommand {\eth} { \end{theorem} }
\newcommand {\ele} { \end{lemma} }

\newcommand {\epr} { \end{proposition} }
\newcommand {\eco} { \end{corollary} }
\newcommand {\ede} { \end{definition} }
\newcommand {\eex} { \end{example} }
\newcommand {\ere} { \end{remark} }

\newcommand \thref[1]{Theorem \ref{t#1}}
\newcommand \leref[1]{Lemma \ref{l#1}}
\newcommand \prref[1]{Proposition \ref{p#1}}

\newcommand \deref[1]{Definition \ref{d#1}}

\newcommand \reref[1]{Remark \ref{r#1}}
\newcommand \lb[1]{\label{#1}}
\def \Rset {{\mathbb R}}         
\def \Cset {{\mathbb C}}

\def \Nset {{\mathbb N}}
\def \g  {\mathfrak{g}}   
\def \ka  {\mathfrak{k}}
  
\def \h  {\mathfrak{h}}
\def \sl {\mathfrak{sl}}
\def \B  {{\mathcal{B}}}           
\def \C  {{\mathcal{C}}}
\def \Ss  {{\mathcal{S}}}
\def \G   {{\mathcal{C}}}
\def \rep {{\bf{\mathrm{rep}}}}
\def \ii   {\bf{1}}
\def \De {\Delta}   
\def \de {\delta}
\def \al {\alpha}
\def \be {\beta}
\def \om {\omega}
\def \e  {\epsilon}   
\def \la {\lambda}
\def \La {\Lambda}
\def \Ps {\Psi}
\def \ph {\varphi}
\def \sig{\sigma}
\def \x  {\xi}
\def \op {\oplus}
\def \mt  {\mapsto}
\def \ra  {\rightarrow}           
\def \lra {\longrightarrow}

\def \hra {\hookrightarrow}
\def \sub {\subset}
\def \st  {\ast}                 
\def \ci  {\circ}

\def \rcor {\rangle}
\def \lcor {\langle}
\def \o  {\otimes}
\def \ol {\overline}            
\def \wt {\widetilde}
\def \wh {\widehat}

\def \D  {{\mathcal{D}}}           
  
\def \id { {\mathrm{id}} }
\def \const { {\mathrm{const}} }

\def \Lie { {\mathrm{Lie}} }
\DeclareMathOperator \Span { {\mathrm{span}} }  
\DeclareMathOperator \ad { {\mathrm{ad}} }  
\DeclareMathOperator \Ker { {\mathrm{Ker}} }
\DeclareMathOperator \tr { {\mathrm{tr}} }
\DeclareMathOperator \qtr { {\mathrm{qtr}} }
\DeclareMathOperator \End { {\mathrm{End}} }
\DeclareMathOperator \Ob { {\mathrm{Ob}} }
\DeclareMathOperator \Hom { {\mathrm{Hom}} }
\DeclareMathOperator \Dom { {\mathrm{Dom}} }
\renewcommand \Im { {\mathrm{Im}} }
\renewcommand \Re { {\mathrm{Re}} }

\begin{document}
\title[Quantum Invariant Measures]
{Quantum Invariant Measures}
\author[N.~Reshetikhin]{Nicolai Reshetikhin}
\thanks{
The first author was partially
supported by NSF grant DMS96-03239}
\address{
Department of Mathematics \\
University of California at Berkeley \\
Berkeley, CA 94720, U.S.A.}
\email{reshetik@math.berkeley.edu}
\author[M.~Yakimov]{Milen Yakimov}
\thanks{
The second author was 
supported by NSF grants
DMS94-00097 and DMS96-03239}
\address{ 
Department of Mathematics \\
University of California at Berkeley \\
Berkeley, CA 94720, U.S.A.}
\email{yakimov@math.berkeley.edu}
\date{}
\begin{abstract}
We derive an explicit expression for the Haar integral on 
the quantized algebra of regular functions $\Cset_q[K]$
on the compact real form $K$ of an arbitrary simply connected
complex simple algebraic group $G.$
This is done in terms of the irreducible $\st$-representations 
of the Hopf $\st$-algebra $\Cset_q[K].$
Quantum analogs of the measures
on the symplectic leaves of the standard Poisson structure on $K$ 
which are (almost) invariant under the dressing action of 
the dual Poisson algebraic group $K^\st$
are also obtained. They are related to the notion of quantum traces
for representations of Hopf algebras. As an application we define
and compute explicitly quantum analogs of Harish-Chandra
$c$-functions associated to the elements of the Weyl group of $G.$
\end{abstract}
\maketitle
\section{Introduction}\lb{Intro}
Let $G$ be a simply connected complex simple algebraic group. 
The cocommutative
Hopf algebra $\Cset[G]$ of regular functions on $G$ has a standard 
quantization, denoted by $\Cset_q[G]$ and called quantized
algebra of regular functions on $G.$ It is a Hopf subalgebra of 
the dual Hopf algebra of the standard quantized universal enveloping
algebra $U_q \g.$ Let $K$ denote a compact real form of $G.$
The complex conjugation in the algebra $\Cset[K](=\Cset[G])$
can be 
deformed to a conjugate linear antiisomorphism $\st$ of $\Cset_q[G]$.
This gives rise to a Hopf $\st$-algebra $(\Cset_q[G], \st)$ 
called quantized algebra of regular functions on $K$ which will
be denoted by $\Cset_q[K].$ 

The Hopf algebra $\Cset_q[K]$ is known \cite{CP} to have 
a unique Haar functional
$H \colon \Cset_q[K] \ra \Cset$ normalized by $H(1)=1.$ 
It is known by a quantum analog of the Schur orthogonality relations.
At the same time an analog of the classical expression for
the bi-invariant functional on $\Cset[K]$ as an integral over $K$ with
respect to the Haar measure was found only in the case of
$SU_2,$ \cite{SV}. The first result which we obtain in this paper
is a representation for the Haar integral on $\Cset_q[K]$
of this type in the general case.

Let us first note that the quantum analog of the set of points 
on $K$ is the set of irreducible $\st$-representations of the 
Hopf $\st$-algebra $\Cset_q[K].$ Its representations were 
classified by Soibelman \cite{Soi} and can be nicely described
by a version of Kirillov--Kostant orbit method.
Fix a maximal torus $T$ of $K.$ Let $G=KAN$ be the related Iwasawa
decomposition of $G.$ The group $K$ has a standard Poisson structure
making it a real Poisson algebraic group which is the semiclassical
structure of the deformation
of $\Cset[K]$ to $\Cset_q[K].$ The double and dual Poisson algebraic 
groups of $K$ are isomorphic to $G$ and $AN$ as real algebraic groups,
respectively.
The dressing action of $AN$ on $K$ is global and is 
explicitly given by the rule \cite{LW, Soi}
\begin{equation}
\lb{dressing}
\de_{a n} (k) \quad \mbox{for} \quad a \in A, n \in N, k \in K 
\quad \mbox{is such that} \quad
a n k= \left( \de_{a n} (k) \right) a_1 n_1
\end{equation}
for some $a_1 \in A, n_1 \in N$ (see \cite{STS, LW} for general
facts about the dressing action). 
Let us choose for each element $w$ of the Weyl group $W$ of $G$
a representative $\dot{w}$ in the normalizer of $A$ in $K.$ 
The orbits of the dressing action of $AN$ on $K$ 
(symplectic leaves of $K)$ 
are $\Ss_{w}. t$ where $w \in W,$ $t \in T$ and $\Ss_{w}$ 
denotes the orbit of $\dot{w}.$ The disjoint union
$\sqcup_{t \in T}\Ss_{w}. t$ is the Bruhat cell $K \cap B \dot{w} B$
where $B$ is the Borel subgroup $B=TAN$ of $G.$
Soibelman proved that the leaves $\Ss_{w}. t$
are deformed to a set $\pi_{w,t}$ of (unequivalent) irreducible 
$\st$-representations of the Hopf $\st$-algebra $\Cset_q[K].$ 
Up to an equivalence they exhaust all such representations of 
$\Cset_q[K].$

Our result on the Haar integral on $\Cset_q[K]$ expresses
it as an integral over the maximal torus $T$ of $K$ of the traces
of the representations $\pi_{w_\ci,t}$ for the maximal
length element $w_\ci$ of $W.$ In other words these are the
irreducible $\st$-representations of $\Cset_q[K]$ 
corresponding to the symplectic leaves in the 
maximal Bruhat cell of $K.$ This result is derived in
Section~5. It is particularly
suited for obtaining integral expressions for quantum spherical functions.
This will be discussed in a future publication.

For each $w \in W$ denote $N_w = N \cap w N_- w^{-1}$ and 
$N_w^+= N \cap w N w^{-1}$ where $N_-$ is the opposite
to $N$ unipotent subgroup of $G.$ Our next result is a quantum analog 
of the Haar measures on the unipotent groups $N_w.$
The symplectic leaf $\Ss_w.t,$ considered 
as an $AN$-homogeneous space 
via the dressing action, is isomorphic to 
\begin{equation}
\lb{s_leaf}
S_w . t = AN/AN_w^+.
\end{equation}
The quotient
$AN/AN_w^+$ does not have a left invariant measure because the 
ratio of the
corresponding modular functions is not equal to 1, see \cite{Hel}. 
Using the factorization $AN= N_w A N_w^+$ we can identify
$AN/A N_w^+ \cong N_w$ which induces a measure on the symplectic leaf
\eqref{s_leaf} from the Haar measure on $N_w.$
The resulting measure transforms under the action of $AN$ 
by the following multiplicative character of $AN$ 
\begin{equation}
\lb{sc_char}
\chi(an) = a^{2(\rho - w \rho)}, \quad a \in A, n \in N.
\end{equation}

The dressing action of $AN =K^\st$ on the symplectic leaf $\Ss_w .t$
of $K$ induces an action of $K^\st$ on the space of functions on 
$\Ss_w .t.$ 
The latter transforms in the quantum situation to an action of 
$\Cset_q[G]$ on the space of linear operators
in the Hilbert space completion $\ol{V}_{w,t}$ 
of the representation space of $\pi_{w,t}.$
It coincides with the standard adjoint action 
\begin{equation}
\lb{adj}
c.L = \sum \pi_{w,t}(c_{(1)}) L \pi_{w,t}(S(c_{(2)})).
\end{equation}
(Here and later we use the standard notation for the comultiplication
in a Hopf algebra 
$\De(c)=\sum c_{(1)} \o c_{(2)}.)$ Let us also note
that
$\Cset_q[K]$ acts by bounded operators in all of 
its $\st$-representations and thus in particular in 
$\ol{V}_{w,t}.$
 
The standard trace in $\ol{V}_{w,t}$ is not a homomorphism from the space of
trace class operators in $\ol{V}_{w,t}$
with the adjoint $\Cset_q[G]$-action \eqref{adj}
to the 1-dimensional representation of $\Cset_q[G]$ determined
by its counit. After Reshetikhin and Turaev such a homomorphism, from
possibly a ``deformation'' of the space of trace class operators,
is called a quantum trace for the Hopf algebra module under 
consideration. We define a space $\B_1^q(\ol{V}_{w,t})$ 
of ``quantum'' trace class operators in $\ol{V}_{w,t},$
stable under the adjoint $\Cset_q[G]$-action \eqref{adj},
and construct a homomorphism from it to the 1-dimensional representation
of $\Cset_q[G]$ determined by a multiplicative character of it which 
is a 
deformation of the character \eqref{sc_char}. Such homomorphisms,
to be called quantum quasi-traces, are treated in Section~6 where 
we also study some of their properties. They are quantum analogs
of the invariant measures on the unipotent groups $N_w$ and
the almost $AN$-invariant measures on the symplectic leaves $\Ss_w. t.$

Section~7 contains an application to quantum analogs of 
Harish-Chandra $c$-functions related to the elements
of the Weyl group of $G.$ They are constructed by the help of 
the quantum quasi-traces from Section~6 and are explicitly computed
by a $q$-analog of the original Harish-Chandra formula.
In the quantum situation 
the role of the factorization formulas for the groups
$N_w$ as products of 1-dimensional unipotent subgroups
is played by tensor product formulas 
for the representations $\pi_{w,t}$
\cite{Soi, KS}.
In a forthcoming publication we will discuss the relation between
the quantum $c$-functions and the asymptotics of quantum spherical 
functions at infinity which is similar to the one in the classical case. 

Sections 2 and 3 review some standard facts about quantized universal 
enveloping algebras, quantized function algebras, and their 
representations. Section~4 deals with 
a family of elements of $\Cset_q[K]$ which enter in  
all formulas for quantum invariant functionals derived in this paper.
\section{Preliminaries on quantized enveloping algebras}\lb{qea}
\subsection{Root data}\lb{rdata}
Let $\g$ be a complex simple Lie algebra of rank $l$ with Cartan matrix
$(a_{ij}).$ Denote by $(.,.)$ the invariant inner product
on $\g$ for which the square length of a minimal root equals 2 in the
resulting identification $\h^\st \cong \h$ for a Cartan subalgebra $\h$
of $\g.$ The sets of simple roots, simple coroots, and fundamental
weights of $\g$ will be denoted by $\{\al_i\}_{i=1}^l,$ 
$\{\al_i\spcheck\}_{i=1}^l,$ and $\{\om_i\}_{i=1}^l,$ respectively.  
Let $P,$ $Q,$ and $Q\spcheck,$ denote the weight, root, 
and coroot lattices of $\g.$ Denote by $\De,$ $\De_+,$ $\De_-,$
and $P_+$ the sets of roots, 
positive/negative roots, and dominant weights of $\g.$
Set $Q_+= \{ \sum m_i \al_i \}$ and 
    $Q_+\spcheck = \{\sum m_i \al_i\spcheck \},$
    $m_i \in \Nset.$ 

Recall that there exists a unique set of relatively prime positive
integers $\{ d_i \}_{i=1}^l$ for which the matrix $(d_i a_{ij})$ is
symmetric and for it 
\[
(\al_i, \al_j) = d_i a_{ij}.
\]

The Weyl group of $\g$ will be denoted by $W.$ The simple reflections
in $W$ will be denoted by $s_i$ and the maximal length element in $W$ by
$w_\ci.$
\subsection{Definition of $U_q \g$}\lb{Uqg-def}
Throughout this paper we will assume that 
$q$ is a real number different from $\pm 1$ and $0.$
The adjoint rational form of the quantized universal enveloping algebra
$U_q \g$ of $\g$ is generated by $K^{\pm 1}_i,$ and $X^\pm_i,$ 
$i=1, \ldots, l,$ subject to the relations  
\begin{gather*}
K_i^{-1} K_i = K_i K^{-1}_i = 1, \, K_i K_j = K_j K_i, 
\\
K_i X^\pm_j K^{-1}_i = q_i^{a_{ij}} X^\pm_j,
\\
X^+_i X^-_j - X^-_j X^+_i = \de_{i,j} \frac{K_i - K^{-1}_i}
{q_i - q^{-1}_i},
\\
\sum_{r=0}^{1-a_{ij}}
\begin{bmatrix} 
1-a_{ij} \\ r
\end{bmatrix}_{q_{i}}
      (X^\pm_i)^r X^\pm_j (X^\pm_i)^{1-a_{ij}-r} = 0, \, i \neq j.
\end{gather*}
It is a Hopf algebra with comultiplication given by
\begin{gather*}
\De(K_i) = K_i \o K_i,
\\
\De(X^+_i)= X^+_i \o K_i + 1 \o X^+_i,
\\
\De(X^-_i)= X^-_i \o 1 + K^{-1}_i \o X^-_i,
\end{gather*}
antipode and counit given by
\begin{gather*}
S(K_i) = K^{-1}_i, \, 
S(X^+_i)= - X^+_i K^{-1}_i, \, 
S(X^-_i)= - K_i X^-_i,
\\
\e(K_i), \, \e(X^\pm_i)=0
\end{gather*}
where $q_i = q^{d_i}.$ As usual $q$-integers, $q$-factorials, and
$q$-binomial coefficients are denoted by 
\[
[n]_q = \frac{q^n - q^{-n}}{q - q^{-1}}, \, 
[n]_q ! = [1]_q \ldots [n]_q, \, 
\begin{bmatrix}
n \\ m  
\end{bmatrix}_{q}
= \frac{[n]_q}{[m]_q [n-m]_q}
\]
for $n, m \in \Nset$ and $m \leq n.$

The conjugate linear antiisomorphism $\st$ of $U_q\g$ defined on its
generators by
\begin{equation}
\lb{st}
K^\st_i = K_i, \,
(X^+_i)^\st = X^-_i K_i, \,
(X^-_i)^\st= K^{-1}_i X^{-}_i
\end{equation}
equips $U_q \g$ with a structure of a Hopf $\st$-algebra.
In the limit $q \ra 1$ the involution $\st$ recovers
the Cartan (anti)involution (conjugate linear antiisomorphism of order 2)
of $\g$ associated to its compact real form $\ka.$
For the definition and properties of Hopf $\st$-algebras
we refer to \cite[pp. 95--97]{KS} and \cite[pp. 117--118]{CP}.

For $i=1, \ldots, l$ the Hopf subalgebra of $U_q \g$ generated by
$K_i$ and $X^\pm_i$ will be denoted by $U_{q_i} \g_i.$ It is naturally
isomorphic to $U_{q_i} \sl_2.$ The canonical embedding
$U_q{\sl_2} \cong U_{q_i} \g_i \hra U_q \g$ will be denoted by $\ph_i.$  
Recall that a $U_q \g$-module is called {\em{integrable}} if the
subalgebras $U_{q_i} \g_i$ act locally finitely.

The subalgebras of $U_q \g$ generated by
$\{K_i\}_{i=1}^l,$ $\{X^+_i\}_{i=1}^l,$ and $\{X^-_i\}_{i=1}^l$ will
be denoted by $U_0,$
$U^+,$ and $U^-,$ respectively. Clearly $U_0$ is a commutative Hopf
subalgebra of $U_q \g$ isomorphic to the group algebra of the
lattice $Q$ equipped with the standard structure of a cocommutative
Hopf algebra.
\subsection{Quantum Weyl group}\lb{braid} 
Let $\B_\g$ denote the (generalized) braid group associated 
to the Coxeter group $W$ with generators $T_i$ corresponding
to the simple reflections $s_i \in W.$ 
For any integrable $U_q \g$-module $V$  
one can define an action of $\B_\g$ on $V.$
It is given by \cite{L}
\[
T_i = \sum_{a, b, c \in \Nset} 
(-1)^b q_i^{ac - b} (X^+_i)^{(a)} (X^-_i)^{(b)} (X^+_i)^{(c)} 
\]
where
\[
(X^\pm_i)^{(n)}= \frac{X^\pm_i}{[n]_{q_i}} \cdot
\]
In the case of the adjoint representation of $U_q \g$ this gives an 
action of the braid group $\B_\g$ on $U_q \g.$ The explicit action of
$T_i$ on the generators $K_j,$ $X^\pm_j$ of $U_q \g$ is
\begin{align*}
&T_i(X^+_i) = - X^-_i K_i, \, 
T_i(X^-_i) = - K^{-1}_i X^+_i, \, 
T_i (K_j) = K_j K^{- a_{ij}}_i, 
\\
&T_i(X_j^+) = \sum_{r=0}^{- a_{ij}} (-1)^r q_i^{-r}
(X_i^+)^{(-a_{ij} -r)} X_j^+ 
(X_i^+)^{(r)} \quad \text{if} \quad i \neq j,
\\ 
&T_i(X_j^-) = \sum_{r=0}^{- a_{ij}} (-1)^r q_i^r
(X_i^-)^{(r)} X_j^-
(X_i^-)^{(- a_{ij}-r)} \quad \text{if} \quad i \neq j.
\end{align*}

The defined actions of $\B_\g$ are compatible in the sense that 
for any integrable $U_q \g$-module $V$ 
\[
T_i . x v = (T_i x) . T_i v, \, \forall x \in U_q \g, \quad
v \in V.
\]

Recall that there exists a canonical section $T \colon W \ra \B_\g$ of the
natural projection $\B_\g \ra W$ (where $T_i \mt s_i).$ If 
\[
w= s_{i_1} \ldots s_{i_n}
\]
is a reduced decomposition of $w \in W$ then the image $T_w$ of
$w$ in $\B_\g$ is defined by
\[
T_w = T_{i_1} \ldots T_{i_n}.
\]
It does not depend on the choice of a reduced decomposition.

The weight subspaces of a $U_0$-module 
(in particular of a $U_q \g$-module) $V$ are defined by
\[
V_\la = \{ v \in V \mid K_i . v = q^{(\la, \al_i)} v \}, 
\quad \la \in P.
\]
The elements of $\B_\g$ preserve the weight space decomposition of an
integrable $U_q \g$-module, in particular
\[
T_w V_\la = V_{w \la}.
\]

\subsection{$R$-matrix}\lb{Rm} Put 
\begin{equation}
\lb{Uk}
U^\pm_k=\bigoplus_{
\begin{smallmatrix}
\la \in \pm Q_+, \\ 
|(\la, \rho\spcheck)| \geq k
\end{smallmatrix}
}
U^\pm_\la, \quad k \in \Nset
\end{equation}
where $\rho\spcheck$ is the half-sum of positive coroots
of $\g.$ Denote by $U_+ \wh{\o} U_-$
the completion of the vector space $U_+ \wh{\o} U_-$ according to
the descending sequence of vector spaces
\[
\left(U^+_k \o U^- \right)\op 
\left(U^+ \o U^-_k \right).
\]
Any element of the completion $U_+ \wh{\o} U_-$ acts 
in the tensor product of two finite dimensional $U_q \g$-modules.

Recall that a representation $V$ of $U_q \g$ is called a type 1
representation if it is a direct sum of its weight subspaces.
For a pair $(V_1,V_2)$ of type 1 $U_q \g$-modules define the
linear operator
$\Ps_{V_1, V_2} \colon V_1 \o V_2 \ra V_1 \o V_2$
by
\[
\Ps_{V_1,V_2}(v_1 \o v_2) = q^{(\la, \mu)} v_1 \o v_2 \quad
\text{if} \quad v_1 \in (V_1)_\la, \, v_2 \in (V_2)_\mu.
\]
Denote also by
$\sig \colon V_1 \o V_2 \ra V_2 \o V_1$ the flip operator
\[
\sig(v_1 \o v_2) = v_2 \o v_1.
\]

There exists \cite{L, KS} a unique element $R \in U^+ \wh{\o} U^-,$ 
called {\em{a quasi $R$-matrix}} for $U_q\g,$ normalized by
\[
R-1 \in U^+_1 \wh{\o} U^-_1 
\]
such that for any pair $(V_1,V_2)$ of finite dimensional 
$U_q \g$-modules of type 1 the composition
\begin{equation}
\lb{quasiR}
\sig \circ \Psi_{V_1, V_2} \circ R: V_1 \o V_2 \ra V_2 \o V_1
\end{equation}
defines an isomorphism of $U_q \g$-modules. 

For any pair $(V_1, V_2)$ of finite dimensional $U_q \g$-modules 
and an element $w \in W$ the actions of $T_w \in \B_\g$ on $V_1, V_2,$
and $V_1 \o V_2,$ to be denoted by $T_{w, V_1}, T_{w, V_2},$ and 
$T_{w, V_1 \o V_2},$ are related as follows. 
There exists a unique element $R^w \in U^+ \wh{\o} U^-$
which does not depend on $V_1$ and $V_2$ 
such that
\begin{equation}
\lb{RwT}
T_{w, V_1 \o V_2} = R^w \left( T_{w, V_1} \o T_{w, V_2} \right).
\end{equation}
As the quasi $R$-matrix $R,$ $R^w$ satisfies
\begin{equation}
\lb{Rw1}
R^w-1 \in U^+_1 \wh{\o} U^-_1.
\end{equation}
The element $R^{w_\ci}$ associated to the maximal element $w_\ci$ of $W$
is equal to the quasi $R$-matrix $R.$
\section{Quantized algebras of functions}
\subsection{Quantized algebras of regular functions} \lb{CqG}
Let $G$ be a
connected, simply connected, complex simple 
algebraic group and $\g = \Lie G.$ The finite dimensional, 
$U_q \g$-modules of type 1 form a quasitensor
category. Hence their matrix coefficients form a Hopf subalgebra of the
Hopf dual $(U_q \g)^\st$ of $U_q \g.$ It is called the quantized algebra
of regular functions on $G$ and is denoted by $\Cset_q[G].$

Every finite dimensional type 1 $U_q \g$-module is a direct sum of
irreducible type 1 $U_q \g$-modules. The latter are highest weight modules
with highest weights $\La \in P_+$ (the corresponding 
module will be denoted by $L(\La)).$ 
The matrix coefficient of $L(\La)$ associated to
$v \in L(\La)$ and $l \in L(\La)^\st$ 
will be denoted by $c^\La_{l, v}:$
\[
c^\La_{l, v} \in \Cset_q[G], \, c^\La_{l, v}(x)= \lcor l, x . v \rcor.
\]
The above implies
\[
\Cset_q[G] = \Span \{ c^\La_{l, v} \mid \La \in P_+, \, v \in L(\La), \, 
l \in L(\La)^\st \}.
\]

The $\st$-involution in $U_q \g$ induces a structure of 
Hopf $\st$-algebra on $\Cset_q[G]$ by
\begin{equation}
\lb{star}
\lcor \xi^\st, x \rcor = \ol{ \lcor \xi, S(x)^\st \rcor },    \quad
\xi \in \Cset_q[G], \, x \in U_q \g.
\end{equation}
The resulting Hopf $\st$-algebra $(\Cset_q[G], \st)$ is called quantized
algebra of regular functions on the compact real form $K$ of $G$
and is denoted by $\Cset_q[K].$ 

The inclusions $\ph_i \colon U_{q_i} \g_i \hra U_q \g$ induce surjective
homomorphisms 
$\ph_i^\st \colon (\Cset_q[G], \st) \ra (\Cset_{q_i}[G_i], \st)$
where $G_i$ is the subgroup of $G$ isomorphic to $SL_2$
with tangent Lie algebra $\g_i$ generated by the root vectors of 
$\pm \al_i.$

We finish this subsection with a simple fact on the explicit
structure of the Hopf $\st$-algebra $\Cset_q[K]$
(see, for instance, \cite[Proposition~13.1.3]{CP}).

Recall that $L(\La)^\st \cong L(-w_\ci \La)$ and 
if we fix these isomorphisms, we can consider
any $v \in L(\La),$ $l \in L(\La)^\st$ as elements of
$L(-w_\ci \La)^\st,$ $L(-w_\ci \La),$ respectively. Recall
that any module $L(\La)$ can be equipped with a unique 
(up to a constant) inner product which turns it into a $(U_q \g, \st)$
$\st$-representation. 

\ble{3.1} (i) The comultiplication, the counit, 
and the antipode of $\Cset_q[G]$ 
are given by
\begin{gather}
\lb{de-cq}
\De(c^\La_{l,v}) = \sum_j c^\La_{l, v_j} \o c^\La_{l_j, v}, \\
\lb{es-cq}
\e(c^\La_{l,v})= \lcor l, v \rcor, \,
S(c^\La_{l,v})= c^{-w_\ci \La}_{v,l}
\end{gather}
where in \eqref{de-cq} $(\{v_j\}, \{l_j\})$ is an 
arbitrary pair of dual bases of $L(\La)$ and $L(\La)^\st.$

(ii) Fix an orthonormal basis $\{v_i\}$ of $L(\La)$ equipped with an
invariant inner product as above and a dual
basis $\{l_j\}$ of $L(\La)^\st.$ The action of the $\st$-involution 
\eqref{star} on the corresponding elements of $\Cset_q[G]$ is given by
\begin{equation}
\lb{st-invol}
(c^\La_{l_i, v_j})^\st= (c^{-w_\ci \La}_{v_i, l_j}).
\end{equation}
\ele
\subsection{Quantized algebra of continuous functions of $K$}Let $G$
be a complex simple algebraic group as in the previous subsection and $K$ 
be its compact
real form. The quantized algebra of continuous functions 
$\C_q(K)$ on $K$ is by definition the $C^\st$-completion 
of the $\st$-algebra $\Cset_q[K]$
with respect to the norm 
\begin{equation}
\lb{sup}
\|f \| = \sup_{\eta} \| \eta(f) \|, \quad f \in \Cset_q[K] 
\end{equation}
where $\eta$ runs through all $\st$-representations of 
$\Cset_q[K].$

The fact that for any $\st$-representation $\eta$ of $\Cset_q[K]$
$\eta(f)$ is a bounded operator and that the supremum in \eqref{sup} is
finite for all $f \in \Cset_q[G]$ follows from the following
identity in $\Cset_q[K]$
\[
\sum_j c^\La_{l_j, v_i} (c^\La_{l_j, v_i})^\st=1
\]
where $\{v_i\}$ and $\{ l_j \}$ are dual bases of $L(\La)$ and
$L(\La)^\st$ as in part (ii) of \leref{3.1}, see
\cite[eq. (13) p. 452]{CP}.

The $C^\st$-algebras $\C_q(K)$ posses natural structures 
of compact matrix quantum groups in the sense of 
Woronowicz \cite{Wo}, see \cite[Section~13.3]{CP}.
\subsection{$\Cset_q[SU_2]$} \lb{sl2}
The $U_q \sl_2$-module $L(\om_1)$
has a basis in which the operators $K_1,$ $X^\pm_1$ act by
\[
K_1 \mt 
\begin{pmatrix}
q & 0 \\
0 & q^{-1}
\end{pmatrix}
, \quad
X^+_1 \mt
\begin{pmatrix}
0 & 1 \\
0 & 0
\end{pmatrix}
, \quad
X^-_1\mt
\begin{pmatrix}
0 & 0 \\
1 & 0
\end{pmatrix}
.
\]  
The corresponding matrix coefficients 
$c_{ij} \in \Cset_q[SL_2]$ $i,j=1,2$ generate 
$\Cset_q[SL_2].$ More precisely:
\ble{Cqsl2}The Hopf algebra $\Cset_q[SL_2]$ is isomorphic to the
algebra generated by $c_{ij},$ $i, j=1,2,$ subject to the relations
\begin{align*}
c_{11} c_{12} = q^{-1} c_{12} c_{11}, & \quad
c_{11} c_{21} = q^{-1} c_{21} c_{11},
\\
c_{12} c_{22} = q^{-1} c_{22} c_{12}, & \quad
c_{21} c_{22} = q^{-1} c_{22} c_{21},        
\\
c_{12} c_{21} = c_{21} c_{12}, & \quad
c_{11} c_{22} - c_{22} c_{11}=
(q^{-1}-q) c_{12} c_{21}        
\\
& \quad
c_{11} c_{22} - q^{-1} c_{12} c_{21}=1.        
\end{align*}
In these generators the comultiplication, the counit, the antipode,
and the $\st$-involution of $\Cset_q[SU_2]$ are given by
\begin{gather*}
\De(c_{ij})= \sum_{k=1,2} c_{ik} \o c_{kj}, \,
\e(c_{ij})= \de_{ij}, 
\\
S(c_{11})=c_{22}, \, 
S(c_{22})=c_{11}, \, 
S(c_{12})=-q c_{12}, \, 
S(c_{21})=-q^{-1} c_{21},  
\\
c_{11}^\st= c_{22}, \, 
c_{21}^\st= -q c_{12}.
\end{gather*}
\ele
A proof of \leref{Cqsl2} can be found, for instance, in 
\cite[Example~2.3.3 and Theorem~3.0.1]{KS}.

Let $q \in \Rset,$ $q > 1.$
The Hopf $\st$-algebra $\Cset_q[SU_2]$ has an infinite dimensional
$\st$-representation $\pi$ on $l^2(\Nset)$ 
given by the following action of
its generators $c_{ij},$ $i,j=1,2$ 
(see \cite{Soi, KS})
\begin{align}
\lb{act1}
&\pi(c_{12}) e_k = q^{-k-1} e_k, \quad 
\pi(c_{11}) e_k = \sqrt{1- q^{-2k}} e_{k-1},
\\
\lb{act2}
&\pi(c_{21}) e_k = -q^{-k} e_k, \quad 
\pi(c_{22}) e_k = \sqrt{1- q^{-2k-2}} e_{k+1}
\end{align}
where $e_{-1}:=0.$
\subsection{Irreducible 
star representations of $\Cset_q[K]$} \lb{reps}
The group of multiplicative characters of the Hopf 
algebra $C_q[G]$ is isomorphic
to the complex torus $(\Cset^\times)^l,$ 
see \cite[Theorem 3.3]{HL} and \cite[Section 10.3.8]{J}
in the case when $q$ is an indeterminate.
The character corresponding
to the $l$-tuple $t =(t_1, \ldots, t_l) \in (\Cset^\times)^l$ 
is given by
\begin{equation}
\lb{char}
\chi_t(c^\La_{l,v})= \prod_{i=1}^l t_i^{(\la, \al_i\spcheck) } 
\lcor l, v \rcor
=\prod_{i=1}^l t_i^{(\la, \al_i\spcheck)} \e(c^\La_{l,v}), \, 
v \in L(\La)_\la.
\end{equation}
The unitary ones among these are the ones corresponding to the real 
torus $(S^1)^l=\{(t_1, \ldots, t_l) \in (\Cset^\times)^l 
\mid |t_i|=1 \}.$ 

  {}  {\em{From now on we will assume that}} $q \in \Rset,$ $q >1.$
Denote by $\pi_i$ the $\st$-representation of 
$(\Cset_{q_i}[G_i], \st) \cong \Cset_{q_i}[SU_2]$ given by 
\eqref{act1}--\eqref{act2}. The $\st$-representation
of $\Cset_q[K] \cong (\Cset_{q}[G], \st)$ induced from it 
by the homomorphism
$\ph^\st_i: (\Cset_{q}[G], \st) \ra (\Cset_{q_i}[G_i], \st)$
will be denoted by $\pi_{s_i}.$ (Recall that $s_i$ denotes the
simple reflection in the Weyl group $W$ of $\g$ corresponding
to the root $\al_i.)$

The irreducible $\st$-representations of the Hopf
$\st$-algebra $\Cset_q[K]$ were classified by Soibelman
\cite{Soi}, see also the book \cite{KS} for an exposition.

\bth{ureps} (i) For any reduced decomposition 
$w=s_{i_1} \ldots s_{i_n}$ of an element $w$ of $W$ 
and any $t \in (S^1)^l$ the tensor product
\begin{equation}
\lb{rep}
\pi_{w, t} = \pi_{s_{i_1}} \o \ldots \o \pi_{s_{i_n}} \o \chi_t
\end{equation}
is an irreducible $\st$-representation of $\Cset_q[K].$

(ii) Up to an equivalence the representation $\pi_{w,t}$ does not 
depend on the choice of reduced decomposition of $w.$

(iii) Every irreducible $\st$-representation of $\Cset_q[G]$
is isomorphic to some $\pi_{w, t}.$
\eth

Denote by $V_{w, t}$ the representation space of $\pi_{w, t}$ equipped
with the Hermitian inner product from \thref{ureps}. The Hilbert space
completion of $V_{w, t}$ with respect to it will be denoted by 
$\ol{V}_{w, t}.$ Then: 

{\em{The representations $\pi_{w, t}$ naturally induce 
irreducible representations of the $C^\st$-algebra $\C_q(K),$
$\pi_{w, t} \colon \C_q(K) \ra \B(\ol{V}_{w, t}).$ The latter 
exhaust all irreducible representations of $\C_q(K)$ up to 
a unitary equivalence.}}

Each module $V_{w, t}$ has a natural orthonormal basis
\begin{equation}
\lb{basis}
e_{k_1, \ldots, k_n} = e_{k_1} \o \ldots \o e_{k_n}\o 1, \quad 
n=l(w), \, k_1, \ldots, k_n \in \Nset 
\end{equation}
induced from the orthonormal basis $\{e_k\}$ of the $\Cset_q[SU_2]$-module
$V$ defined by \eqref{act1}--\eqref{act2}. Here $1$ denotes a 
(fixed) vector of the 1-dimensional representation of 
$\Cset_q[G]$ corresponding to $\chi_t.$

For an element $w$ of the Weyl group $W$ denote by $I_w$ the $\st$-ideal
of $\Cset_q[K]$ generated by
\begin{equation}
\lb{I_w}
c^\La_{l, v_\La} \quad \text{such that} \quad
\La \in P_+, \, \lcor l, U^+ T_w . v_\La \rcor = 0
\end{equation}
where $v_\La$ denotes a highest weight vector of $L(\La).$

The annihilation ideals of the representations $\pi_{w, t}$ are contained
in $I_w$ \cite{Soi,KS}:
\begin{equation}
\lb{annih}
\ker \pi_{w, t} \sub I_w.
\end{equation}
\section{A family of elements $a_{\La,w} \in \Cset_q[K]$}
\subsection{Definitions} For a dominant integral weight $\La \in P_+$
and a highest weight vector $v_\La$ of $L(\La)$
denote by $l_{\La, w}$ the unique element of $L(\La)^\st_{-w \La}$ 
such that
\[
\lcor l_{\La, w}, T_w v_\La \rcor =1.
\]
(The uniqueness follows from the fact that $\dim L(\La)_{w \La}=1.)$ 

Define
\begin{equation}
\lb{defa}
a_{\La, w} = c^\La_{l_{\La, w}, v_\La}.
\end{equation}
Note that $a_{\La, w}$ does not depend on the choice of highest
weight vector $v_\La$ of $L(\La).$

The $\st$-subalgebras of $\Cset_q[K]$ generated by $a_{\La, w}$ played an
important role in Soibelman's classification 
of the irreducible $\st$-representations of $\Cset_q[K],$ 
see \thref{ureps}.  
Most of the results in this subsection are due to 
Soibelman \cite{Soi}.
We include their proofs since \cite{Soi} does not assume 
the normalization made in the definition of $a_{\La, w}.$ 

Properties \eqref{RwT} and \eqref{Rw1} of the elements 
$R^w \in U^+ \wh{\o} U^-$ allow to write $l_{\La, w}$ and 
thus $a_{\La, w}$ slightly more explicitly.
Let $l_\La=l_{\La, 1},$ i.e. let 
$l_\La \in L(\La)^\st_{-\La}$ be the unique element such that
\[
\lcor l_\La, v_\La \rcor =1.
\] 
Then \eqref{RwT}, \eqref{Rw1} imply 
\[
l_{\La, w}=T_w l_\La
\]
and thus
\begin{equation}
\lb{defa2}
a_{\La, w} = c^\La_{T_w l_\La, v_\La}.
\end{equation}

\bpr{aalg} (i) The elements $a_{\La, w},$ $a^\st_{\La, w} \in \Cset_q[K],$
$\La \in P_+$ are normal modulo $I_w:$
\begin{align}
\lb{norm}
a_{\La, w} c^{\La'}_{l, v} - q^{(\La, \la')-(w\La, \mu')}
c^{\La'}_{l, v} a_{\La, w} \in I_w, 
\\
\lb{norm*}
a^\st_{\La, w} c^{\La'}_{l, v} - q^{(\La, \la')-(w\La, \mu')}
c^{\La'}_{l, v} a^\st_{\La, w} \in I_w, 
\end{align}
for $v \in L(\La')_{\la'},$ $l \in L(\La')^\st_{-\mu'}.$

(ii) The images of $\{a_{\La, w}, a^\st_{\La, w} \}_{\La \in P_+}$
in $\Cset_q[K]/I_w$ generate a commutative subalgebra.
More precisely the following identity holds in $\Cset_q[K]$
\begin{equation}
\lb{proda}
a_{\La_1, w} a_{\La_2, w}= a_{\La_1+\La_2, w}, \quad \forall
\La_1, \La_2 \in P_+.
\end{equation}
\epr

Proofs of \prref{aalg} can be found in \cite{Soi, KS}.
The property \eqref{norm} follows from the existence of a 
quasi $R$-matrix for $U_q \g,$ see \eqref{quasiR}.
Eq. \eqref{norm*} follows from \eqref{norm}, \leref{3.1}, and
the fact that the ideals $I_w$ are stable under the $\st$-involution.
The first statement in part (ii) is a direct consequence of
part (i). The second statement in (ii) follows from 
the existence of the element $R^w \in U^+_1 \wh{\o} U^-_1$ with the
properties \eqref{RwT}, \eqref{Rw1}
and the fact that $v_{\La_1} \o v_{\La_2} \in L(\La_1) \o L(\La_2)$
generates a submodule isomorphic to $L(\La_1+\La_2).$
\subsection{The action of $a_{\La, w}$ in $V_{w,t}$}
\ble{comA}Let $w, w' \in W$ be such that
$w=s_i w'$ and $l(w) = l(w') +1$
for some simple reflection $s_i \in W.$ Then 
\[
\De( a_{\La, w}) - 
c^\La_{l_{\La, w}, T_{w'} v_\La} \o a_{\La, w'}.
\]
\ele
\begin{proof} According to \eqref{de-cq} $\De(a_{\La, w})$ is given by
\[
\De( a_{\La, w}) =
\sum_j c^\La_{l_{\La, w}, v_j} \o c^\La_{l_j, v_\La}
\]
where $(\{v_j\}, \{l_j\})$ is a pair of dual bases of $L(\La)$ and 
$L(\La)^\st$ consisting of weight vectors $(v_j \in L(\La)_{\la_j},$ 
$l_j \in L(\La)_{-\la_j},$ $\la_j \in P).$
The definition \eqref{I_w} of $I_{w'}$ implies
\begin{equation}
\lb{1w}
c^\La_{l_j, v_\La} \in I_{w'} \quad
\text{if} \quad \la_j \notin w \La + Q_+.
\end{equation}
The map $\ph^\st_i \colon \Cset_q[G] \ra \Cset_{q_i}[G_i]$ acts on 
the matrix coefficients of a $U_q \g$-module 
by restricting the module to $U_{q_i} \g_i.$
Since $w=s_i w'$ and $l(w)=l(w')+1$
\[
w^{-1} \al_i\spcheck \in - Q_+\spcheck.
\]
If $\La$ is a dominant weight, then 
\[
\lcor \La, w^{-1} \al_i\spcheck \rcor \leq 0 \quad \text{and thus} 
\quad 
\lcor w \La, \al_i\spcheck \rcor \leq 0.
\]
Hence $T_w v_\La$ is a lowest weight vector for 
the $U_{q_i} \g_i$-submodule of $L(\La)$ 
generated by $T_w v_\La.$ The corresponding 
$U_{q_i} \g_i$-highest weight vector
is $T_{w'} v_\La$ and
\begin{equation}
\lb{2w}
c^\La_{l_j, v_\La} \in (\ph^\st_i)^{-1}(T_{s_i}) \quad
\text{if} \quad \la_j \notin \{ w \La, w \La + \al_i, \ldots, w' \La \}.
\end{equation}
The lemma now follows from \eqref{1w} and \eqref{2w}.
\end{proof}

For an element $w \in W$ and a reduced decomposition 
$w = s_{i_1} \ldots s_{i_n}$ of it denote
\begin{equation}
\lb{w_j}
w_j= s_{i_{j+1}} \ldots s_{i_n}, \, j=0, \ldots, n-1, \quad
w_n=1.
\end{equation}

\bpr{ra} In the notation \eqref{w_j} the action of the elements 
$a_{\La, w}$ in the module $V_{w, t}$ is given by
\begin{equation}
\lb{aV1}
\pi_{w, t}(a_{\La, w})=
\bigotimes_{j=1}^n
\pi_{s_{i_j}}(a_{(w_j \La, \al_{i_j}\spcheck) \om_{i_j},
s_{i_j}}) . \prod_{i=1}^l t_i^{(\La, \al_i\spcheck)}.
\end{equation}
In the orthonormal basis $\{e_{{k_1}, \ldots, {k_n}} \}_{k_j=0}^\infty,$ 
of $V_{w, t},$ see \eqref{basis}, the elements $a_{\La, w}$ act 
diagonally by 
\begin{equation}
\lb{aV2}
\pi_{w, t}(a_{\La, w}) . e_{{k_1}, \ldots, {k_n}}
= \prod_{j=1}^n q^{-(k_j+1)(w_j \La, \al_{i_j})}
\prod_{i=1}^l t_i^{(\La, \al_i\spcheck)}
 e_{{k_1}, \ldots, {k_n}}.
\end{equation}
\epr

Formula \eqref{aV1} follows by induction from \leref{comA} and the
definition \eqref{char} of the multiplicative characters $\chi_t$ of 
$\Cset_q[G].$ 
To prove \eqref{aV2} we first compute that in $\Cset_q[SL_2]$
\begin{equation}
\lb{relat}
a_{\om_1, s_1}= - q c_{21}
\end{equation}
(cf. Section~\ref{sl2})
and then use \eqref{proda} which implies 
$a_{m \om_1, s_1}= (a_{\om_1, s_1})^m.$ We also use the identity
$d_i \al_i\spcheck = (\al_i, \al_i) \al_i\spcheck/2 = \al_i,$
see Section~\ref{rdata}. 
\section{The Haar integral on $\C_q(K)$}
\subsection{Definition and the Schur orthogonality relations}
Recall that a left invariant integral on a Hopf algebra $A$ is a
linear functional $H$ on $A$ satisfying 
\begin{equation}
\lb{left_Haar}
(\id \o H) \left( \De(a) \right) = H(a). 1, \quad
\forall a \in A.
\end{equation}
Analogously is defined a right invariant integral. 
In the analytic setting a left Haar integral for
a $C^\st$-Hopf algebra $A$ is a state $H$ on $A$ satisfying
\eqref{left_Haar}, see \cite{Wo}.

\bpr{uHaar}There exists a unique left invariant integral $H$ on the
Hopf algebra $\Cset_q[K]$ normalized by $H(1)=1.$ It is also 
right invariant and can be uniquely extended to a bi-invariant 
Haar integral on $\C_q(K).$ It is given by a quantum version of 
the classical Schur orthogonality relations:
\[
H( c^\La_{l, v} c^{ \La'}_{l', v'})=
\frac{ \de_{\La, \La'} \lcor l, v' \rcor \lcor l', v \rcor}
{\sum_\la \dim L(\La)_\la q^{2 (\la, \rho)} }   
\]
or equivalently by
\begin{equation}
\lb{Schur}
H(c^\La_{l, v})= \de_{\La, 0} \lcor l, v \rcor.
\end{equation}
\epr
\subsection{Statement of the main result}
\bth{Haar}The bi-invariant integral $H$ on $\C_q(K)$ 
$(q \in \Rset, \; q >1)$ is given
in terms of the irreducible representations $\pi_{w, t}$ of 
$\C_q(K)$ by
\begin{equation}
\lb{haa}
H(c) = \left( \prod_{\be \in \De_+}
(q^{(2 \rho, \be)}-1) \right)
\int_{(S^1)^l}
\tr_{\ol{V}_{w_\ci, t}} 
( \pi_{w_\ci, t} ( a_{\rho, w_\ci} a^\st_{\rho, w_\ci} c )) 
d t
\end{equation}
where $w_\ci$ is the maximal length element of the Weyl group $W$ 
of $\g,$ $\rho$ is the half sum of all positive roots of $\g,$
and $d t$ is the invariant measure on the torus $(S^1)^l$ 
normalized by $\int_{(S^1)^l} d t =1.$ 
\eth

In the special case of $K=SU_2$ \thref{Haar} was established by
Soibelman and Vaksman
\cite{SV}. Similar formula is also
known for quantum spheres \cite{Sinf, SVsov}. 
\thref{Haar} answers Question 3 in \cite{Sinf}. 

Note that formula \eqref{aV2} implies that 
$\pi_{w, t}(a_{\rho, w_\ci} a^\st_{\rho, w_\ci})$
is a trace class operator in $\ol{V}_{w, t}.$ Since 
$\pi_{w, t}(c)$ is a bounded operator 
for $c \in \C_q(K),$ the product is also a trace class
operator in $\ol{V}_{w, t}$ for all $c \in \C_q(K).$ From the definition
\eqref{rep} of $\pi_{w, t}$ 
it is also clear that 
\[
\tr_{\ol{V}_{w_\ci, t}} ( \pi_{w_\ci, t} 
(a_{\rho, w_\ci} a^\st_{\rho, w_\ci} c) )
\]
is a continuous function in $t \in (S^1)^l$
for a fixed $c \in \C_q(K)$ and that the r.h.s. of
\eqref{haa} defines a continuous linear functional on $\C_q(K).$

By identifying 
$(S^1)^l \cong \{ (t_1, \ldots, t_l) \in \Cset^l \colon |t_i| =1 \},$
the normalized invariant measure on the torus $(S^1)^l$ is 
represented as
\[
dt = \frac{1}{(2 \pi i)^l}
\frac{ dt_1}{t_1} \wedge \ldots \wedge \frac{d t_l}{t_l} \cdot
\]

In Sections~\ref{reduce} and \ref{sl2case} we show that the functional
$\wt{H}$ on $\Cset_q[G]$ given by the right hand side of \eqref{haa}
satisfies
\begin{equation}
\lb{0}
\wt{H}(c^\La_{l, v})= 0 \quad \text{if} \quad \La \neq 0.
\end{equation}
In Section~\ref{constant} we check that it satisfies the normalization
condition $\wt{H}(1)=1.$ Combined with \eqref{Schur} this proves
\thref{Haar}.
\subsection{Proof of \eqref{0}: reduction to the rank 1 case} \lb{reduce}
Recall first the following simple characterization of $w_\ci \in W.$
\ble{wci}The maximal length element $w_\ci \in W$ is the only element $w 
\in W$
that has a representation of the form $w = w' s_i$ with $l(w')= l(w)-1$
for an arbitrary simple reflection $s_i.$
\ele

\leref{wci} follows from the so called ``deletion condition'',
see \cite{H}, and the property of $w_\ci$ that it is the only element
$w \in W$ such that $w^{-1}(\al_i)$ is a negative root of $\g$ for all
simple roots $\al_i$ of $\g.$

We show that \eqref{0} for $K=SU_2$ implies its validity in the general
case. Let $\La \in P_+,$ $\La \neq 0.$ Equip $L(\La)$ with a Hermitian
inner product making it a $(U_q \g , \st)$ $\st$-representation, recall
\eqref{st}. 
Denote
\[
L_i =\{ v \in L(\La) \mid U_{q_i} \g_i . v = 0 \}, \, i=1, \ldots, l.
\]
Since $L(\La)$ is an irreducible $U_q \g$-module
\[
\cap_{i=1}^l L_i = 0
\]
and
\[
L_1^\perp + \ldots + L_2^\perp 
=(\cap_{i=1}^l L_i)^\perp = L(\La).
\]
Hence to show \eqref{0} it is sufficient to show that
\begin{equation}
\lb{reduce0}
\wt{H}(c^\La_{l, v}) =0 \quad \text{if} \quad
v \in L_m^\perp \quad \text{for some} \quad m = 1, \ldots, l.
\end{equation}
Note that $L_m^\perp$ is the span of the nontrivial irreducible 
$U_{q_m}\g_m$-submodules of $L(\La).$

Choose a reduced decomposition of $w_\ci$ of the form
\[
w_\ci = s_{i_1} \ldots s_{i_{n_\ci-1}} s_m
\]
and consider the corresponding model for the representation
$\pi_{w_{n_\ci}, t}$ 
\[
\pi_{w_{n_\ci}, t} \cong \pi_{s_{i_1}} \o \ldots \o \pi_{i_{n_\ci-1}}
\o \pi_{s_m} \o \chi_t.
\]
Taking trace over the component $\pi_{s_m} \o \chi_t$ of 
$\pi_{w_{n_\ci}, t}$ and using 
\eqref{de-cq} and \eqref{aV1} we see that to prove
\eqref{reduce0} it is sufficient to prove that 
\begin{equation}
\lb{intbr}
\int_{(S^1)^l}\tr_{\ol{V}_{s_m}}
(\pi_{s_m} (a_{\om_m, s_m} a^\st_{\om_m, s_m} \ph^\st_m(c^\La_{l', v}) )) 
d t =0
\quad \text{for all} \quad
l' \in L(\La)^\st.
\end{equation}
(Recall that by definition $(w_\ci)_{n_\ci}=1,$ see \eqref{w_j}.)
Since $v \in L^\perp_m$
\[
\ph^\st_m(c^\La_{l, v}) =
\sum_{p} c^{p \om_m}_{l_p, v_p}
\]
with all $p > 0.$ By appropriately breaking the integral 
\eqref{intbr} into a product of a 1-dimensional and an $(l-1)$-dimensional
integrals
one sees that \eqref{intbr} follows from \eqref{0} for $K=SU_2.$ 
\subsection{Proof of \eqref{0}: the case of $\Cset_q[SU_2]$} \lb{sl2case}
Our proof in the rank 1 case is similar to the one from \cite{SV}.
\leref{Cqsl2} implies that $\Cset_q[SU_2]$ is spanned by the 
elements
\[
c_{11}^m c_{12}^p c_{21}^r \, \text{and} \, 
c_{22}^m c_{12}^p c_{21}^r \quad \text{for} \quad m, p, r \in \Nset.
\]
The Haar functional $H$ acts on them by \cite[Example 13.3.9]{CP}
\[
H(c_{11}^m c_{12}^p c_{21}^r)=
H(c_{22}^m c_{12}^p c_{21}^r)=
\de_{m, 0} \de_{p, r} \frac{(-q)^p(q^2-1)}{q^{2p+2}-1} \cdot
\]
We check that the functional $\wt{H}$ has the same property. 
This implies \eqref{0} for $K=SU_2.$

Recall from \eqref{relat} that $a_{\om_1, s_1} = -q^{-1} c_{21}$
and thus $a_{\om_1, s_1}^\st = c_{12},$ see \leref{Cqsl2}. 
Using \eqref{act1}--\eqref{act2} we compute
\begin{align*}
\tr_{\ol{V}}(\pi(a_{\om_1, s_1} 
a_{\om_1, s_1}^\st c_{ii}^m c_{12}^p c_{21}^r)) 
&= \de_{m, 0} \sum_{k=0}^\infty - q^{-1}.q^{-(k+1)(p+1)}.(-q^{-k})^{r+1} 
\\
&= \de_{m, 0} \frac{(-q)^r}{q^{p+r+2}-1}
\end{align*}
for $i=1,2.$ This gives
\begin{align*}
&\frac{1}{2\pi i}
\int_{S^1} \tr_{\ol{V}_{s_1, t}}
(\pi_{s_1, t}( a_{\om_1, s_1} a_{\om_1, s_1}^\st 
  c_{ii}^m c_{12}^p c_{21}^r)) \frac{d t}{t} 
\\
&= \de_{m, 0} \frac{(-q)^r}{q^{p+r+2}-1} 
\int t^{r-p-1} dt =
\de_{m, 0} \de_{p, r} \frac{(-q)^p}{q^{2p+2}-1}
\end{align*}
$(i=1,2)$ which shows that $\wt{H} =H$ in the case $K=SU_2.$
\subsection{Checking the normalization $\wt{H}(1)=1$} \lb{constant}
Let $w_\ci=s_{i_1} \ldots s_{i_{n_\ci}}$ be a reduced decomposition
of the maximal element of $W.$ Using \eqref{aV2} and the notation
\eqref{w_j} we compute
\begin{align*}
\int_{(S^1)^l}
\tr_{\ol{V}_{w_\ci, t}} (\pi_{w_\ci, t}
(a_{\rho, w_\ci} a^\st_{\rho, w_\ci}) ) 
d t &=\prod_{j=1}^{n_\ci}
\left(
\sum_{k_j=0}^\infty
q_{i_j}^{-2(k_j+1)((w_\ci)_j \rho, \al_{i_j}\spcheck)}
\right)
\\
&=\prod_{j=1}^{n_\ci} 
\frac{1}{1-q_{i_j}^{-(2\rho, (w_\ci)^{-1}_j \al_{i_j}\spcheck)}} \cdot
\end{align*}
Note that $q_{i}^{(\la, \al_i\spcheck)}=q^{(\la, \al_i)}$
for all simple roots $\al_i$ of $\g.$
The set of elements $(w_\ci)^{-1}_j \al_{i_j} \in Q,$
$j=1, \ldots, n_\ci,$
coincides with the set of positive roots of $\g.$ This together with
the definition of the functional $\wt{H}$ via the r.h.s. of \eqref{haa}
gives 
\[
\wt{H}(1)=1.
\]
\subsection{Semiclassical limit} \lb{intsc}
Here we explain the semiclassical 
analog of the integral formula from \thref{Haar}.

As earlier $G$ denotes a complex simple algebraic group 
and $K$ denotes a compact real form of $G.$
For each element $w$ of the Weyl group $W$ of $K$ choose 
a representative $\dot{w}$ of it in the normalizer of a fixed maximal
torus $T$ of $K.$ Using the related Iwasawa decomposition of 
$G,$ introduce the map
\begin{equation}
\lb{a_w}
a_w \colon N \ra A \quad \mbox{by} \quad \dot{w}^{-1} n \dot{w}
= k_1 a_w(n) n_1, \quad k_1 \in K, n_1 \in N,
\end{equation}
see for instance \cite{Lu}.
It can be pushed down to a well defined map from the symplectic
leaf $\Ss_w$ to $A$ 
\[
a_w(\de_n \dot{w}) := a_w(n), \; n \in N.
\] 
We refer to the introduction for details
on the dressing action of $AN$ on $K$ related to the standard 
Poisson structure on $K.$

The semiclassical analog of formula \eqref{haa} is the 
following formula for the normalized Haar integral on $K$ 
\begin{equation}
\lb{semicl_int}
H(f) = \left( \prod_{\be \in \De_+} \frac{\pi}{(\rho, \be)} \right)
\int_{\Ss_{w_\ci} \times T} a_{w_\ci}(k)^{-2\rho} f(k. t) \mu_{w_\ci} dt, 
\quad f \in C(K).
\end{equation}
Here $\mu_{w_\ci}$ denotes the Liouville volume form on the 
symplectic leaf 
$\Ss_{w_\ci}$ corresponding to the maximal element $w_\ci \in W$ and
$d t$ denotes the invariant measure on the torus $T$
normalized by $\int_T dt =1.$ Recall that $\Ss_{w_\ci} \times T$ 
coincides with the maximal Bruhat cell of $K.$ 

Formula \eqref{semicl_int} can be easily proved following the
idea of Sections \ref{reduce}--\ref{constant}
on the basis of the product formulas \cite{Soi, KS} for 
the symplectic leaves $\Ss_w$ of $K,$ $w \in W$
\begin{equation}
\lb{prod-sl}
\Ss_w= \Ss_{s_{i_1}} \ldots \Ss_{s_{i_n}}
\end{equation}
where $s_{i_1} \ldots s_{i_n}$ is a reduced decomposition of $w.$

The integral with respect to the symplectic measure on
the leaf $\Ss_{w}. t$ is (up to a factor) a semiclassical limit of 
the trace in the module $\ol{V}_{w,t}.$ 

At the end we explain the connection between the functions
$a_w^{-2\rho}$ on the leaves $\Ss_w$ and the operators
$\pi_{w, t}(a_{\rho, w} a_{\rho, w}^\st)$ in $\ol{V}_{w,t}.$
Let us consider the highest weight module $L(\La)$ of
$\g$ with heighest weight $\La$ and the matrix coefficient
\[
a_{\La, w} \in \Cset[G], \quad
a_{\La, w}(g) = \lcor l_{w, \La}, g v_\La \rcor, 
\, g \in G
\]
where $v_\La$ is a highest weight vector of $L(\La)$ and
$l_{w, \La}\in L(\La)^\st_{-w \La}$ is normalized by 
$\lcor l_{w, \La}, \dot{w} v_\La \rcor = 1,$ cf. \eqref{defa}.
It is easy to show that the restriction of 
$a_{\La, w}$ to the symplectic leaf $\Ss_w$ coincides 
with $a_w^{-\La}$
\[
a_{\La, w}|_{\Ss_w}= a_w^{-\La}.
\]
For $t \in T$ the functions 
\[
|a_{w, \rho}(k.t)|^2=|a_{w, \rho}(k)|^2= a_w(k)^{-2\rho}, 
\; k \in \Ss_w
\]
are semiclassical analogs of the linear operators
$\pi_{w, t}(a_{w, \rho} a_{w, \rho}^\st)$ in $\ol{V}_{w,t}.$ 
\section{Quantum quasi-traces of $V_{w,t}$}
\subsection{Motivation} Let $A$ be a Hopf algebra
and $A^\st$ be its dual Hopf algebra.
Denote by $A^\ci$ the dual Hopf algebra of $A$  
equipped with the opposite comultiplication. 
Recall \cite{Dr, CP} that the quantum double
$\D(A)$ of $A$ is isomorphic to $A \o A^\ci$ as a coalgebra
and the following commutation relation holds in $\D(A)$
\begin{equation}
\x a = \sum \lcor \x_{(1)}, a_{(3)} \rcor \x_{(2)} a_{(2)}
       \lcor S^{-1} \x_{(3)}, a_{(1)} \rcor, 
\quad \x \in A^\st, a \in A.
\lb{double}
\end{equation}
Analogously to the classical situation one defines
a quantum dressing action $\de$ of $A^\st$ on $A.$ Using the
identification $\D(A) \cong A \o A^{\st}$ as vector
spaces, set
\[
\de_\x a= (\id \o \e)( \x a).
\] 
In view of the commutation relation \eqref{double} it is explicitly
given by
\[
\de_\x a = \sum \lcor \x_{(1)}, a_{(3)} \rcor a_{(2)}
       \lcor S^{-1}\x_{(2)}, a_{(1)} \rcor.
\]
It is dual to the standard adjoint action of $A^\st$ on itself
\[
\ad_\x \x' = \sum \x_{(1)} \x' S(\x_{(2)})
\]
in the sence that
\begin{equation}
\lb{adA}
\lcor \ad_\x \x', a \rcor =
\lcor \x', \de_{S (\x)} a \rcor.
\end{equation}

For any representation $\pi$ of $A^\st$ in the vector space 
$V$ the adjoint action of $A^\st$ on itself lifts to an
action of $A^\st$ in the space of linear operators on $V$
by
\begin{equation}
\lb{mod_ad}
\ad_\x L = \sum \pi(\x_{(1)}) L \pi(S \x_{(2)}).
\end{equation}

Suppose that $A^\st$ is a deformation of the Poisson Hopf algebra
$\Cset[F]$ of regular functions on a Poisson algebraic group $F.$
According to Kirillov--Kostant orbit method philosophy 
an irreducible $A^\st$-module $V$ can be viewed as a quantization
of a symplectic leaf $\Ss$ in $F.$ The left action of $A^\st$ in the space
of linear operators in $V$ is a deformation
of the Poisson $\Cset[F]$-module of functions on the leaf $\Ss.$
At the same time the dual Poisson algebraic group $F^\st$ of
$F$ acts in the space of functions on $\Ss$ by the dressing action.
The quantum analog of this action is the 
the adjoint action \eqref{mod_ad} of $A^\st$ in the space of
linear operators in the $A^\st$-module $V.$ This leads to:

{\em{The quantum analog of a measure on the symplectic leaf $\Ss$ 
in the Poisson algebraic group $F$
which is invariant up to a multiplicative character of $F^\st$ 
is a homomorphism from a subspace of linear operators in the $A^\st$-module 
$V,$ equipped with the $A^\st$-action \eqref{mod_ad}, to a 1-dimensional 
representation of $A^\st.$}} 

In the next subsection we will develop this idea from a
categorical point of view and relate it to the notion of quantum traces
for $A^\st$-modules. In analogy, the defined more general morphisms
will be called quantum quasi-traces. Subsections 6.3 and 6.4 construct
such morphisms for the irreducible $\st$-representations of
the quantized algebras of functions $(\Cset_q[G], \st).$

\subsection{Definitions}
Let $\G$ be a $\Cset$-linear, rigid, monoidal category with 
identity object $\ii.$ 
Recall that $\G$ is called balanced if for each object
$V \in \Ob(\G)$ there exists an isomorphism
\[
b_V \colon V \stackrel{\cong}{\ra} V^{\st \st}
\]
such that 
\begin{gather}
\lb{bal1}
b_{V_1} \o b_{V_2}= b_{V_1 \o V_2}, \\ 
\lb{bal2}
b_{V^\st} =(b_V^\st)^{-1}, \\
\lb{bal3}
b_{\ii} = \id_{\ii}.
\end{gather}
Given a Hopf algebra $C$ over the field $\Cset$ let $\rep_C$ denote 
the category of its finite dimensional modules equipped with 
the left dual object $V^\st$ of $V \in \Ob(\G)$ defined by
\begin{equation}
\lb{dua}
\lcor c . \x, v \rcor = \lcor \x, S(a) . v \rcor, 
\quad \x \in V^\st, v \in V.
\end{equation}
The spaces $\Hom_{\Cset}(V_1, V_2),$ $V_1, V_2 \in \Ob(\G)$ can be
equipped with the canonical $C$-action 
\begin{equation}
\lb{ad}
c . L = \sum \pi_{V_1}(c_{(1)}) L \pi_{V_2}(S(c_{(2)})), 
\quad L \in \Hom_{\Cset}(V_1, V_2).
\end{equation}
Here the Hopf algebra $C$ plays the role of the Hopf algebra
$A^\st$ from the motivation in the previous subsection, cf. 
\eqref{mod_ad} and its derivation from the quantum dressing action. 

Clearly 
\[
\Hom_{\Cset}(V_1, V_2) \cong V_2 \o V_1^\st
\]
as $C$-modules. In particular, for this action 
$\Hom_{\Cset}(V, \e)$
is canonically isomorphic to $V^\st$ where, by abuse of notation,
$\e$ denotes the 1-dimensional representation of 
$C$ defined by its counit.

Reshetikhin and Turaev \cite{RT} defined the following notion 
of quantum trace for a finite dimensional $C$-module $V.$ 

\bde{qtr}A {\em{quantum trace}} for a finite dimensional $C$-module $V$ 
is a homomorphism
\[
\qtr_V \colon \End_{\Cset}(V) \ra \e
\]
of $C$-modules for the action of $C$ on $\End_{\Cset}(V)$ defined in
\eqref{ad}.
\ede

The pairing
\[
\End_\Cset(V) \cong V \o V^\st \ra \Cset
\]
is not a homomorphism of $C$-modules where $\Cset$
is given the structure of the $C$-module corresponding to 
the counit $\e.$ At the same time the opposite pairing 
\[
V^\st \o V \ra \Cset
\]
has this property.
If $\rep_C$ is balanced each $V \in \Ob(\G)$ has a quantum trace
defined by the composition \cite{RT}
\[
\End_{\Cset}(V) \cong 
V \o V^\st \stackrel{b_V \o \id}{\longrightarrow} 
V^{\st \st} \o V^\st \ra \e
\]
or explicitly
\[
\qtr_V(L) = \tr_V(b_V L), \, L \in \End(V).
\]  
Here $b_V$ is considered as a linear endomorphism of $V$
using the canonical identification of $V$ and $V^{\st \st}$
as vector spaces. 

The properties \eqref{bal1}--\eqref{bal2} of the
balancing morphisms $b_V$ imply the following properties
of the quantum traces $\qtr_V$
\begin{gather}
\lb{qtr1}
\qtr_{V_1 \o V_2}(L_1 \o L_2) =
\qtr_{V_1}(L_1) \qtr_{V_2}(L_2), \\
\lb{qtr2}
\qtr_{V^\st}(L^\st) = \qtr_{V}(L) 
\end{gather}
for all $L_i \in \End_\Cset(V_i).$

In \cite{RT, RT2} it was proved that the category of finite dimensional
type 1 $U_q \g$-modules is balanced and this was used for constructing 
invariants of links and 3-dimensional manifolds.

We would like to incorporate in \deref{qtr} the possibility
for an invariant up to a character ``quantum measure'', 
as explained in the previous
section, and the general case of an infinite dimensional $C$-module 
$V.$ We will restrict ourselves to representations of $C$
$\pi \colon C \ra \B(V)$ by bounded operators in a Hilbert
space $V$ and will call them {\em{bounded}} representations of $C.$ 
The Hermitian inner product in $V$ is not assumed to posses 
any invariance properties and the linear operators $\pi(c),$
$c \in C,$ in $V$ are not assumed to be uniformly bounded. 
The dual $V^\st$ of 
such a bounded representation $\pi \colon C \ra \B(V)$
is defined in the Hilbert space $V^\st$ of bounded functionals 
on $V$ by formula \eqref{dua}. Obviously it is again a bounded 
representation.

\bde{wequiv}
Two bounded representations of a Hopf algebra $C$ 
$\pi_i \colon C \to \B(V_i)$ in the Hilbert spaces 
$V_i,$ $i=1,2,$ will be called weakly equivalent if 
$V_i$ contain dense $C$-stable subspaces $W_i \sub V_i$ 
which are equivalent as $C$-modules.
\ede

The point here is that the equivalence can be given by an unbounded
operator $b \colon W_1 \stackrel{\cong}\ra W_2$ which therefore
does not extend to the full space $V_1.$

\bde{balanced} A bounded representation $\pi \colon C \to \B(V)$ 
of a Hopf algebra $C$ in a Hilbert space $V$ will be called 
quasi-balanced if there exists a multiplicative character $\chi$ 
of $C$  for which $V$ and $\chi \o V^{\st \st}$ are weakly equivalent.
\ede
By abuse of notation we denote by $\chi$ the 1-dimensional
$C$-module corresponding to the multiplicative character $\chi$
of $C.$

In other words the bounded $C$-module $V$ is balanced if there exists 
an invertible linear operator $b_V$ in $V$ with dense domain and range 
such that $\Dom b_V$ is $C$-stable and  
\begin{equation}
\lb{bal}
b_V \pi(c) = \sum 
\chi(c_{(1)}) \pi(S^2(c_{(2)}))b_V, \quad \forall c \in C.
\end{equation}
(Here we use the canonical identification of $V^{\st \st}$ and $V$
as Hilbert spaces.)

\bre{restdual}
Often $V$ is the Hilbert space completion of a $C$-module $W,$
equipped with a Hermitian inner product, which is a direct sum of 
mutually orthogonal finite dimensional submodules $W_\mu$
for a Hopf subalgebra $B$ of $C$
\begin{equation}
\lb{dec}
W= \op_{\mu} W_\mu.
\end{equation}
The restricted dual of such a module $W$ with respect to the
decomposition \eqref{dec} as a direct sum of finite dimensional 
subspaces is naturally a $C$-module of the same type.
The double restricted dual $W^{\st \st}$ of $W$ 
is canonically isomorphic to $W$ as a vector space.

If $W \cong \chi \o W^{\st \st}$ as $C$-modules then the modules
$V$ and $\chi \o V^{\st \st}$ are weakly equivalent and $V$ is
a quasi-balanced $C$-module.
\ere

Let $\pi \colon C \to \B(V)$ be a quasi-balanced representation 
as above. We call the subspace of the space of linear operators in $V$
with dense domains
\[
\B^q_1(V) := \B_1(V) b^{-1}_V
\]
a space of {\em{quantum trace class operators}} in the $C$-module $V.$
Here $\B_1(V)$ stands for the standard trace class in $V.$ 
It is naturally a $C$-module by
\[
c. L = \sum \pi(c_{(1)}) L \pi(S( a_{(c)}))
\]
because $C$ acts in $V$ by bounded operators. The linear map
$\qtr_V \colon \B^q_1(V) \ra \Cset$ given by
\[
\qtr_V(L) := \tr_V( L b_V )
\]
is a well defined homomorphism of $C$-modules
\[
\qtr_V \colon \B^q_1(V) \ra \chi.
\]
It will be called a quantum {\em{quasi-trace}} for the module $V.$ 

\bre{quasi2}One can as well use the space 
\[
\wt{\B}^q_1(V) := b^{-1}_V \B_1(V) 
\]
instead of $\B^q_1(V).$ When $b^{-1}_V$ is not defined on
the full space $V$ the composition $b^{-1}_V L_0,$ $L_0 \in \B_1(V)$ 
need not have a dense domain in $V.$ Because of this, 
it is convenient to use the space $\wt{\B}^q_1(V)$ only
when $b_V^{-1}$ has full domain. In that case the space 
$\wt{\B}^q_1(V)$ is
also a $C$-module and the following map
\[
\qtr_V \colon \wt{\B}^q_1(V) \ra \chi, \,
\qtr_V(L) := \tr( b_V L )
\]
is a homomorphism of $C$-modules. 
\ere

\bre{balanalg} It is natural to look for a quasi-balancing map 
$b_V \in \End_\Cset(V)$ for a bounded representation 
$\pi \colon C \ra \B(V)$ of the form
\[
b_V= \pi(a_V)
\]
for some $a_V \in C.$ The definition \eqref{bal} implies that such a map
$\pi(a_V)$ provides a quasi-balancing endomorphism if 
$\pi(a_V)$ is an invertible linear operator in $V$ with a dense
range satisfying
\begin{equation}
\lb{balalg}
a_V c - \sum \chi(c_{(1)}) S^2(c_{(2)})a_V \in \Ker \pi,
\quad \forall c \in C
\end{equation}
for some multiplicative character $\chi$ of $A.$

Thus quasi-balancing of the modules of a Hopf algebra $A$ is related to
the properties of the the square of the antipode $S$ of $A.$ This is
analogous to the usual case of balancing when $\chi=\e$ and \eqref{balalg}
reduces to
\[
a_V c = S^2(c)a_V \in \Ker \pi, \quad \forall c \in C,
\]
see \cite{RT}.

Similarly 
\[
b_V = \pi_V(a_V)^{-1}
\]
is a quasi-balancing map for the $C$-module $V$ if
$\pi_V(a_V)$ is an invertible operator in $V$ with a dense range
such that
\begin{equation}
\lb{balalg-1}
c a_V -  \sum  
\chi(c_{(1)}) a_V S^2(c_{(2)}) \in \Ker \pi, \quad
\forall c \in C.
\end{equation}
\ere
\subsection{Main construction}
In this subsection we construct quasi-balancing morphisms
for the $\Cset_q[G]$-modules $V_{w, t}.$
As  was pointed out in Section~\ref{reps} they are bounded 
$\Cset_q[G]$-modules in the terminology from the previous subsection.

Set
\[
2\rho= \sum_{\al \in \De_+} \al= \sum_{i=1}^l p_i \al_i
\]
for some positive integers $p_i$ and denote
\begin{equation}
\lb{rho}
q^{2 \rho} = \prod_{i=1}^l K_i^{p_i} \in U_q \g.
\end{equation}
Its commutation with the generators $X^\pm_i$ of
$U_q \g$ is given by
\[
q^{2 \rho} X_i^\pm q^{-2 \rho} = q^{\pm(2 \rho, \al_i)} X_i^\pm, \quad 
\forall i=1, \ldots, l.
\]

As it is well known the square of the antipode in $U_q \g$ is given by
the following lemma. 
\ble{S2Uq} For all $x \in U_q\g$ 
\[
S^2(x) = q^{2 \rho} x q^{- 2 \rho}.
\]
\ele

For an arbitrary element 
$\nu= \sum_i m_i \al_i\spcheck$ of the coroot
lattice $Q\spcheck$ of $\g$
we set
\begin{equation}
\lb{qnu}
q^\nu := (q^{m_1}, \ldots, q^{m_l}) \in (\Cset^\times)^l
\end{equation}
and consider the multiplicative character 
$\chi_{q^\nu}$ of $\Cset_q[G].$ 
It is explicitly given by
\[
\chi_{q^\nu}(c_{l, v}^\La) =
q^{(\nu, \mu)} \lcor l, v \rcor, \quad l \in L(\La)^\st_{-\mu},
\]
recall \eqref{char}.

  {}  From \leref{S2Uq} we deduce the following
properties of $S^2$ in $\Cset_q[G].$

\ble{S2Cq} {\em{(}}i{\em{\/)}} If $v \in L(\La)_\la$ 
and $l \in L(\La)^\st_{-\mu}$ then
\begin{equation}
\lb{S2}
S^2( c^\La_{l, v}) = q^{2(\rho, \la-\mu)} c^\La_{l, v}.
\end{equation}

{\em{(}}ii{\em{\/)}} For all elements $w \in W$ 
\[
c a_{\rho, w} a_{\rho, w}^\st - \sum 
\chi_{q^{2(w\rho-\rho)}}(c_{(1)}) a_{\rho, w} a^\st_{\rho, w} 
S^2(c_{(2)}) \in I_w, \quad \forall c \in \Cset_q[G],
\]
recall \eqref{I_w}.
\ele
\begin{proof} (i) By a straightforward computation, for all $x \in U_q \g$
\[
\lcor S^2(c^\La_{l, v}), x \rcor =
\lcor c^\La_{l, v}, S^2(x) \rcor =
\lcor c^\La_{l, v}, q^{2 \rho} x q^{-2 \rho} \rcor
= q^{2(\rho, \la-\mu)} \lcor c^\La_{l, v}, x \rcor.
\]

(ii) Combining part (i) with the identities \eqref{norm} and \eqref{norm*}
gives
\[
c^\La_{l, v} a_{\rho, w} a_{\rho, w}^\st-
q^{2(w \rho- \rho, \mu)} a_{\rho, w} a_{\rho, w}^\st S^2(c^\La_{l, v})
\in I_w, 
\quad \forall l \in L(\La)_{-\mu}^\st, v \in L(\La). 
\]
which implies \eqref{S2} in view of \eqref{de-cq}.
\end{proof}

Let us fix an element $w \in W,$ a reduced decomposition
$w=s_{i_1} \ldots s_{i_n}$ of it, and an element $t \in (S^1)^l.$ 
Consider the $(\Cset_q[G], \st)$-module $V_{w,t}.$ 
We will make use of the 
notation \eqref{w_j}
\[
w_j = s_{i_{j+1}} \ldots s_{i_n}, \, j=0, \ldots, n-1,
\quad w_n =1  
\]
and of the basis $e_{k_1, \ldots, k_n},$ $k_j \in \Nset$
of $V_{w,t}$ from \eqref{basis}.

Formula \eqref{aV2} implies that the space $V_{w,t}$ decomposes
as a sum of weight subspaces with respect to the action
of the commutative subalgebra of $\Cset_q[G]$ spanned by
$a_{\La, w},$ $\La \in P_+$ (recall part (ii) of \prref{aalg}) as 
\begin{equation}
\lb{wps}
V_{w,t} = \bigoplus_{\mu \in Q_+}
\Span \{ e_{k_1, \ldots, k_n} \mid 
\sum_{j=1}^n (k_j+1) w_j^{-1} \al_{i_j}= \mu \}.
\end{equation}
All weight subspaces of $V_{w,t}$ are finite dimensional 
and we can identify the corresponding double restricted 
dual $V_{w,t}^{\st \st}$ with $V_{w,t}$ as a vector space. 

Part (ii) of \leref{S2Cq} and the fact that the ideal
$I_w$ contains the annihilation ideal of $V_{w,t},$ see
\eqref{annih}, imply that
$\pi_{w,t}(a_{\rho, w} a_{\rho, w}^\st)^{-1}
\colon V_{w,t} \ra V_{w,t}$ induces an isomorphism 
of the $\Cset_q[G]$-modules $V_{w,t}$ and $\chi \o V_{w,t}^{\st \st}.$
In view of \reref{restdual}, 
$b_{w,t}=\pi_{w,t}(a_{\rho, w} a_{\rho, w}^\st)^{-1}$
defines a quasi-balancing map for the $\Cset_q[G]$-module $\ol{V}_{w,t}.$
Explicitly in the basis \eqref{basis} of $V_{w,t},$
$\pi_{w,t}(a_{\rho, w} a_{\rho, w}^\st)^{-1}$ acts
diagonally by
\begin{equation}
\lb{nor_act}
\pi_{w,t}(a_{\rho, w} a_{\rho, w}^\st)^{-1}.
e_{k_1, \ldots, k_n}=
\prod_{j=1}^n q^{2(k_j+1)(w_j \rho, \al_{i_j})}
e_{k_1, \ldots, k_n},
\end{equation}
recall \eqref{aV2}.

Define the set of quantum trace class operators in
the $\Cset_q[G]$-module $\ol{V}_{w,t}$ by
\begin{equation}
\lb{BqV}
\B^q_1(\ol{V}_{w,t})=
\B_1(\ol{V}_{w,t}) 
\pi_{w,t}(a_{\rho, w} a_{\rho, w}^\st).
\end{equation} 
It is clear from \eqref{nor_act} that 
$\pi_{w,t}(a_{\rho, w} a_{\rho, w}^\st)$ is a compact 
operator and thus
\[
\B^q_1(\ol{V}_{w,t}) \sub \B_1(\ol{V}_{w,t}).
\]
Using \prref{aalg}, observe that
\begin{equation}
\lb{squad}
\pi_{w,t}(a_{2\rho, w} a_{2\rho, w}^\st) 
=\pi_{w,t}(a_{\rho, w} a_{\rho, w}^\st)^2.
\end{equation}

Finally define the quantum quasi-trace functional
$\qtr_{\ol{V}_{w,t}} \colon \B^q_1(\ol{V}_{w,t}) \ra \Cset$ by
\begin{equation}
\lb{trV}
\qtr_{\ol{V}_{w,t}} (L) = 
\const_w
\tr_{\ol{V}_{w,t}}( L \pi_{w,t}(a_{\rho, w} a_{\rho, w}^\st)^{-1})
\end{equation}
where
\begin{equation}
\lb{const}
\const_w = \prod_{\be \in \De_+ \cap w^{-1} \De_-}
(q^{(2 \rho, \be)}-1).
\end{equation}

\bpr{qtrVwt} The $\Cset_q[G]$-modules $\ol{V}_{w,t}$ are quasi-balanced
with multiplicative characters 
$\chi_{2(w \rho - \rho)}$ and quasi-balancing morphisms
$b_{w,t}=\pi_{w,t}(a_{\rho, w} a_{\rho, w}^\st)^{-1}.$ 
The space of quantum trace class operators in $\ol{V}_{w,t}$
and quantum quasi-trace morphisms
\[
\qtr_{\ol{V}_{w,t}} \colon \B^q_1(\ol{V}_{w,t}) \ra 
\chi_{2(w \rho - \rho)}
\]
are given by \eqref{BqV} and \eqref{trV}. The morphisms
$\qtr_{\ol{V}_{w,t}}$ are normalized by 
\begin{equation}
\lb{normtr}
\qtr_{\ol{V}_{w,t}} (\pi_{w,t}(a_{2\rho, w} a_{2\rho, w}^\st) )=1.
\end{equation}
\epr

To check \eqref{normtr} it is sufficient to check that 
\[
\tr_{\ol{V}_{w,t}} (\pi_{w,t}(a_{\rho, w} a_{\rho, w}^\st)) =
\prod_{\be \in \De_+ \cap w^{-1} \De_-}
(q^{(2\rho, \be)}-1)^{-1},
\]
recall \eqref{squad}. 
This easily follows from \eqref{nor_act} using the standard fact
\begin{equation}
\lb{roots}
\{ w_j^{-1} \al_{i_j} \}_{j=1}^n =
\De_+ \cap w^{-1} \De_-
\end{equation}
in the notation of \eqref{w_j}, see for instance \cite{H}.

\bre{plmeasure} Consider again the compact group $K$ equipped with 
the standard Poisson structure, see the introduction and
Section~\ref{intsc}. Recall the notation $N_w = N \cap w N_- w^{-1}$ and 
$N_w^+ = N \cap w N w^{-1},$ $w \in W,$ where $N_-$ is the unipotent
subgroup of $G$ which is dual to $N$ with respect to the fixed
complex torus $TA$ of $G.$ The symplectic leaf $\Ss_w. t$ of $K,$
considered as an $AN$ homogeneous space under the dressing action, 
is isomorphic to $AN /A N_w^+.$ 
We choose as a base point of $\Ss_w. t$ the point $\dot{w}. t.$ 

Denote by $\mu_{w, t}$ the Liouville volume form on the leaf $\Ss_w. t.$
The diffeomorphisms
\begin{equation}
\lb{SwNw}
\Ss_w. t \cong AN /A N_w^+ \cong N_w
\end{equation}
induce a measure $d n_w$ on $\Ss_w. t$ from the Haar measure on 
$N_w.$ The second one comes from the factorization
$AN = N_w A N_w^+.$ The measure $d n_w$ will be normalized by 
\[
\int \Big{|} a_{w, 2 \rho}|_{\Ss_{w}. t} \Big{|}^2 d n_w=1,
\]
cf. Section~\ref{intsc}.

The relation between the volume forms $\mu_{w,t}$ and $d n_w$ on
$\Ss_w. t$ was found by Lu \cite{Lu}. It is given by
\begin{equation}
\lb{Lu_rel}
d n_w = \prod_{\be \in \De_+ \cap w^{-1} \De_-} \left(
\frac{(\rho, \be)}{\pi} \right)
\Big{|} a_{w, \rho}|_{\Ss_{w}. t} \Big{|}^{-2} \mu_{w,t}.
\end{equation}
It is easy to compute that the measure $d n_w$ on $\Ss_w. t$ 
transforms under the dressing action 
of $AN=K^\st$ by
\[
\de_{an} \left( d n_w \right) = a^{2(\rho - w \rho)} d n_w.
\]

The quantum quasi-trace morphisms
\[
\qtr_{\ol{V}_{w,t}} \colon \B^q_1(\ol{V}_{w,t}) \ra 
\chi_{2(w \rho - \rho)}
\]
are quantum analogs of the measures $d n_w$ on $\Ss_w. t$
and thus also of the Haar measures on the unipotent subgroups
$N_w$ of $G.$ The traces in the modules $V_{w,t}$ can be considered 
as quantizations of the Liouville volume forms $\mu_{w, t}$ on the leaves
$\Ss_w. t.$ The relation \eqref{trV} is a quantum version of Lu's
relation \eqref{Lu_rel}.
\ere
\subsection{Tensor product properties of the quasi-balancing morphisms
$b_{w,t}$} 
\lb{mult}
When $w, w' \in W$ are such that $l(w w')= l(w) + l(w')$ the tensor
product of $(\Cset_q[G], \st)$-modules $V_{w, t} \o V_{w', t'}$ is again
an irreducible $(\Cset_q[G], \st)$-module, see \leref{piwt} 
below. Here we discuss the
relation between the corresponding quasi-balancing morphisms
constructed in the previous subsection.

For an element $t=(t_1, \ldots, t_l) \in (\Cset^\times)^l$ denote
its $j$-th component by $(t)_j:=t_j.$
Define an action of the Weyl group $W$ of $\g$ on the torus 
$(\Cset^\times)^l$ by
\[
(w(t))_i := \prod_j t_j^{m_{ij}} \quad
\text{where} \quad w^{-1} \al_j\spcheck = \sum_i m_{ij} 
\al_i\spcheck.
\] 
It can be easily identified with the conjugation action of $W$ on a 
complex torus of $G.$
It is straightforward to check that 
\[
\chi_{w(t)}(c^\La_{l,v})= \prod_{i=1}^l t_i^{(\la, w^{-1}\al_i\spcheck) }
\lcor l, v \rcor,
\]
cf. \eqref{char}.

Fix $w \in W$ and a reduced decomposition 
$w = s_{i_1} \ldots s_{i_n}$ of it. The representation space 
$V_{w, t},$ recall \thref{ureps}, 
is canonically identified with the vector space
\[
V_w = V_{s_1} \o \ldots \o V_{s_n}
\] 
for all $t \in (\Cset^\times)^l.$ 
(As earlier we will not show explicitly the dependence on the choice of a
reduced decomposition of $w.)$
Under this identification
the basis \eqref{basis} of $V_{w, t}$ corresponds to the basis
\begin{equation}
\lb{basis_vw}
e_{k_1, \ldots, k_n} = e_{k_1} \o \ldots \o e_{k_n}, \quad
n=l(w), \, k_1, \ldots, k_n \in \Nset
\end{equation}
of $V_w.$

In the notation \eqref{w_j} define the linear operator $J_{w, t}$
in $V_w$ acting diagonally in the above basis of $V_w$ by
\begin{equation}
\lb{Joper}
J_{w, t} . e_{k_1, \ldots, k_n} =
\prod_{j=1}^n ( w_{j-1}(t) w_j(t^{-1}) )_{i_j}^{k_j+1}
e_{k_1, \ldots, k_n}.
\end{equation}

\ble{inter} For all $w \in W,$ and $t, t' \in (\Cset^\times)^l$ the
operator $J_{w, t'}$ defines an isomorphism of the 
$\Cset_q[G]$-representations $\chi_{w(t')} \o \pi_{w,t}$
and $\pi_{w,t} \o \chi_{t'} \cong \pi_{w,t t'}$ in the natural
identification of their representation spaces with $V_w.$
\ele

\leref{inter} is checked directly in the case of $G = SL_2$
using the defining identities \eqref{act1}--\eqref{act2} 
for the $\Cset_q[SL_2]$-module $\pi,$ see Section \ref{sl2}. 
This implies the lemma when $w$ is a
simple reflection and the general case is proved by induction on $l(w).$

\ble{piwt} Let $w, w' \in W$ be such that $l(ww') = l(w) +l(w')$
and $t, t' \in (S^1)^l.$ The linear operator $J_{w', (w')^{-1}(t)}$
induces the unitary equivalence of $(\Cset_q[G], \st)$-modules
\begin{equation}
\lb{Pi}
\Pi_{w, t; w', t'} \colon \ol{V}_{w, t} \o \ol{V}_{w', t'} \ra 
\ol{V}_{w w', (w')^{-1}(t) t'}
\end{equation}
by identifying the spaces 
$V_{w,t} \o V_{w',t'} \cong V_w \o V_{w'} \cong V_{ww'}
\cong V_{w w', (w')^{-1}(t) t'}.$ (The product of two reduced
decompositions of $w$ and $w'$ is used as a reduced decomposition 
of $w w'.$)
\ele

In the setting of \leref{piwt} the $\Cset_q[G]$-module
$\ol{V}_{w w', (w')^{-1}(t) t'}$ admits a quasi-balancing
morphism constructed from the quasi-balancing morphisms
$b_{w,t}$ and $b_{w',t'}$ for the modules $\ol{V}_{w, t}$ and
$\ol{V}_{w', t'}.$ It is given by the composition  
\begin{align}
\lb{comp}
&V_{w w', (w')^{-1}(t) t'}
\stackrel{\Pi_{w, t; w', t'}^{-1}}{\lra}
V_{w, t} \o V_{w', t'}  
\stackrel{\id \o b_{w', t'}}{\lra}
\\
\nn
&V_{w, t} \o \chi_{q^{2(w'\rho-\rho)}} \o V_{w', t'}^{\st \st}  
\stackrel{J^{-1}_{w, q^{2(w'\rho-\rho)}} \o \id}{\lra}
\chi_{q^{2w(w'\rho-\rho)}} \o V_{w, t} \o V_{w', t'}^{\st \st}  
\stackrel{b_{w, t} \o \id}{\lra}
\\
\nn
&\chi_{q^{2w(w'\rho-\rho)}} \o \chi_{q^{2(w\rho-\rho)}} \o 
V_{w, t}^{\st \st} \o V_{w', t'}^{\st \st}  
\stackrel{\Pi_{w, t; w', t'}^{\st \st}}{\lra}
\chi_{q^{2(ww'\rho-\rho)}} \o V_{w w', (w')^{-1}(t) t'}^{\st \st}.
\end{align}
The restricted duals of the modules $V_{w,t}$ and $V_{w',t'}$
are taken with respect to the weight space decomposition \eqref{wps}
for the commutative subalgebras of $\Cset_q[G]$ spanned
by $a_{\La, w}$ and $a_{\La, w'},$ $\La \in P_+,$ respectively.
Recall also the notation \eqref{qnu}.

\bpr{equal} If $w, w' \in W$ are such that $l(ww') = l(w) +l(w')$
and $t, t' \in (S^1)^l$ then the quasi-balancing map for the 
$\Cset_q[G]$-module $\ol{V}_{w w', (w')^{-1}(t) t'}$ given by the
composition \eqref{comp} coincides with the quasi-balancing
map $b_{w w', (w')^{-1}(t) t'}.$
\epr

To prove \prref{equal} observe that in the natural identification
of the representation spaces in \eqref{comp} with
$V_w \o V_{w'}$ the composition is simply
\[
b_{w,t} J^{-1}_{w, q^{2(w'\rho-\rho)}} \o b_{w', t'}.
\]
(We use again the product of two reduced decompositions of $w$ and $w'$
as a reduced decomposition of $w w'.)$ 
Now the proposition easily follow from \eqref{nor_act} and the 
following formula for the action of $J^{-1}_{w, q^{2(w'\rho-\rho)}}$
in the basis \eqref{basis_vw} of $V_{w,t}$ which is a 
direct consequence from \eqref{Joper}
\[
J^{-1}_{w, q^{2(w'\rho-\rho)}}. e_{k_1, \ldots, k_n}
=\prod_{j=1}^{l(w)}q^{2 (k_j+1) (w_j (w' \rho -\rho), \al_{i_j})}
e_{k_1, \ldots, k_n}.
\]

This computation implies also the following connection between the spaces
$\B_1^q(\ol{V}_{w,t}),$ $\B_1^q(\ol{V}_{w',t'})$ and 
$\B_1^q(\ol{V}_{w w', (w')^{-1}(t) t'})$ when $l(w w') = l(w) + l(w').$

\bco{corol} If $L \in \B_1^q(\ol{V}_{w,t})$ and
$L' \in \B_1^q(\ol{V}_{w',t'})$ in the setting of \prref{equal},
then 
\[
LJ_{w, q^{2(w'\rho-\rho)}} \o L'
\in \B_1^q(\ol{V}_{w w', (w')^{-1}(t) t'})
\]
and
\[
\qtr_{\ol{V}_{w w', (w')^{-1}(t) t'}}
(LJ_{w, q^{2(w'\rho-\rho)}} \o L')=
\frac{\const_{w w'}}{\const_w \const_{w'} } 
\qtr_{\ol{V}_{w, t}}(L) \qtr_{\ol{V}_{w', t'}}(L')
\]
where $\const_w$ is given by \eqref{const}. 
\eco
\section{An application: quantum Harish-Chandra $c$-functions}
Denote by 1 the identity element $(1, \ldots, 1)$ of the real torus
$(S^1)^l.$ 
According to \eqref{aV2} the linear operators 
$\pi_{w,1}(a_{\om_i, w})$ in $\ol{V}_{w,1}$ are compact, selfadjoint
with spectrum contained in $[0, \infty).$ For different values of $i$ they
mutually commute.

Hence for each $\la \in \h$ we can 
define the linear operator in $\ol{V}_{w,1}$
\[
d_{\la, w} = \prod_{i=1}^l \pi_{w,1}(a_{\om_i, w})^{\la_i}
\]
where $\la_i = (\la, \al_i\spcheck),$ i.e. $\la = \sum \la_i \om_i.$
It is obvious that 
\begin{equation}
\lb{eq1}
d_{\la, w}= \pi_{w,1}(a_{\la, w})
\quad \text{when} \quad \la \in P_+ \sub \h
\end{equation}
and
\begin{equation}
\lb{eq2}
d_{\la_1, w} d_{\la_2, w}= 
d_{\la_1+\la_2, w}, \quad
\forall \la_1, \la_2 \in \h.
\end{equation}

\ble{dqrt}The linear operator $d_{i \la + 2 \rho, w}$ in $\ol{V}_{w,1}$
is quantum trace class (belongs to $\B^q_1(\ol{V}_{w,1}))$ 
if and only if 
\[
\Im (\la, \be) < 0, \quad \forall \be \in \De_+ \cap w^{-1}\De_-.
\]
\ele

\begin{proof} The operator $d_{i \la + 2 \rho, w}$ in 
$\ol{V}_{w,1}$ belongs to $\B^q_1(\ol{V}_{w,1})$ if and only if
\[
d_{i \la, w} \in \B_1(\ol{V}_{w,1})
\]
because of \eqref{eq1}, \eqref{eq2}, and the  
the selfadjointness of $\pi_{w,1}(a_{\rho, w}).$
The operator $d_{i \la, w}$ is diagonal in the orthonormal basis 
\eqref{basis} of $\ol{V}_{w,1}$ and according to
\eqref{aV2} acts by
\begin{equation}
\lb{diag_il}
d_{i \la, w}.
e_{k_1, \ldots, k_n}=
\prod_{j=1}^n q^{-i(k_j+1)(w_j \la, \al_{i_j})}
e_{k_1, \ldots, k_n},
\end{equation}
recall the notation \eqref{w_j}. It is clear that
the linear operator $d_{i\la, w}$ in $\ol{V}_{w,1}$ is trace class
if and only if $\Re (i\la, w^{-1}_j \al_{i_j}) > 0$ for
$i=1, \ldots, n=l(w)$ which implies the statement 
because of \eqref{roots}.
\end{proof}

\bde{qcHCh}The function 
\begin{equation}
\lb{quant_c}
c_{w^{-1}}^q(\la) = \qtr_{\ol{V}_{w,1}}(d_{i\la+2 \rho, w}) =
\tr_{\ol{V}_{w,1}}(d_{i\la, w})
\end{equation}
in the domain 
$\{ \la \in \h \mid \Im (\la, \be) <0, \: \forall 
\be \in \De_+ \cap  w^{-1} \De_-\}$
will be called quantum Harish-Chandra $c$-function associated to the 
element $w^{-1}$ of the Weyl group $W$ of $\g.$
\ede

\bpr{comp}For all $w \in W$ the quantum Harish-Chandra $c$-function
$c_w^q(\la)$ is given by
\[
c_w^q(\la) =
\prod_{\be \in \De_+ \cap w \De_-}
\frac{q^{(2 \rho, \be)}-1}{q^{(i \la, \be)} -1} \cdot
\]
\epr

This proposition follows from \eqref{diag_il} and \eqref{roots}
similarly to the proof of the normalization \eqref{normtr}.

\bre{semicl_cfunct} \prref{comp} is a quantum analog of 
the Harish-Chandra formula for the $c$-function in the 
case of complex simple Lie groups, generalized later by 
Gindikin and Karpelevich to arbitrary real reductive groups. 

Recall the setting of Section~\ref{intsc}
and \reref{plmeasure}. Let $d n_w$ denote the Haar measure 
on the unipotent subgroup $N_w$ of $G.$ 
The classical Harish-Chandra $c$-function
associated to the element $w^{-1} \in W$ is given 
by the integral formula
\[
c_{w^{-1}}(\la) = \int_{N_w} a_w(n)^{-(i\la+2 \rho)} d n_w, 
\quad \la \in \h, \: \Im (\la, \be) <0, \, \forall \be \in 
\De_+ \cap w^{-1}\De_-,
\]
recall the definition \eqref{a_w} of the map $a_w \colon N \ra A.$
We refer to \cite{Hel} for a detailed treatment of spherical functions
and to \cite{Lu} for an interpretation of the $c$-function
in terms of the Poisson geometry of $K,$ see in particular Example~2.8 in 
\cite{Lu}.

The linear operators $d_{i\la+2\rho}$ in the modules $\ol{V}_{w,1}$ 
can be thought of as quantizations of the pushforwards of the
functions $a_w(n)^{-(i\la+2 \rho)}$ on $N_w$ to the symplectic leaves
$\Ss_w$ by the dressing action, using the base 
points $\dot{w} \in \Ss_w$
(i.e. using the diffeomorphisms \eqref{SwNw}). As was explained in 
\reref{plmeasure} the quantum quasi-traces 
$\qtr_{\ol{V}_{w, 1}}$ in the $\Cset_q[G]$-modules $\ol{V}_{w,1}$ are
quantizations  of the pushforwards of the Haar measures
on $N_w$ to the symplectic leaves $\Ss_w.$

The classical Harish-Chandra formula 
\[
c_w(\la) =
\prod_{\be \in \De_+ \cap w \De_-}
\frac{(2 \rho, \be)}{(i \la, \be)}, \quad
\la \in \h, \: \Im (\la, \be) <0, \, \forall \be \in 
\De_+ \cap w^{-1}\De_-
\]
is proved by induction on the length of $w,$ see
\cite[Chapter IV, \S 6]{Hel}. Lu \cite{Lu} found that this argument
is essentially based on the product formula 
\eqref{prod-sl} for the leaves $\Ss_w.$ Our computation 
relies on its quantum counterpart -- the tensor product formula \eqref{rep}
for the representations $\pi_{w,t},$ cf. also Section~\ref{mult}. 
\ere
    
\end{document}